# Geometrically Intrinsic Nonlinear Recursive Filters I: Algorithms

R. W. R. Darling[1] [2]

[1] P. O. Box 535, Annapolis Junction, Maryland 20701-0535, USA. E-mail: rwrd@afterlife.ncsc.mil

[2] Supported by the Air Force Office of Scientific Research

**ABSTRACT:** The Geometrically Intrinsic Nonlinear Recursive Filter, or GI Filter, is designed to estimate an arbitrary continuous-time Markov diffusion process $X$ subject to nonlinear discrete-time observations. The GI Filter is fundamentally different from the much-used Extended Kalman Filter (EKF), and its second-order variants, even in the simplest nonlinear case, in that:

- It uses a quadratic function of a vector observation to update the state, instead of the linear function used by the EKF.

- It is based on deeper geometric principles, which make the GI Filter coördinate-invariant. This implies, for example, that if a linear system were subjected to a nonlinear transformation $f$ of the state-space and analyzed using the GI Filter, the resulting state estimates and conditional variances would be the push-forward under $f$ of the Kalman Filter estimates for the untransformed system - a property which is not shared by the EKF or its second-order variants.

The noise covariance of $X$ and the observation covariance themselves induce geometries on state space and observation space, respectively, and associated canonical connections. A sequel to this paper develops stochastic differential geometry results – based on "intrinsic location parameters", a notion derived from the heat flow of harmonic mappings – from which we derive the coördinate-free filter update formula. The present article presents the algorithm with reference to a specific example – the problem of tracking and intercepting a target, using sensors based on a moving missile. Computational experiments show that, when the observation function is highly nonlinear, there exist choices of the noise parameters at which the GI Filter significantly outperforms the EKF.

**RUNNING HEAD: :** Geometrically Intrinsic Nonlinear Filters

# 1 Background on Nonlinear Filtering

## 1.1 Example: A Nonlinear Filtering Problem in Target Tracking

D'Souza, McClure, and Cloutier [8], [9] consider the following tactical air-to-air missile intercept problem. The state of the target is represented by a position, velocity, and acceleration in space, making nine dimensions in all (the authors also model three time constants as state variables). Data consists of a sequence of noisy observations of: range, angle from vertical, azimuth, and range-rate, all measured from a missile with known position, velocity, and acceleration. The goal of filtering in this case is to provide a sequence of "good" estimates of the state of the target, based on all measurements so far, so as to defeat the target's possible evasive maneuvers and intercept it.

D'Souza et al [8] point out that, although the state dynamics can be modeled linearly, the observations are a highly nonlinear function of the state (see Section 2.7.a). Alternatively, if a spherical coördinate frame, based on the missile, is used, then observations are linear, but the state dynamics are highly nonlinear. Moreover the Extended Kalman Filter, and standard second-order filters, will give a different set of answers in the Cartesian coördinate frame than in the spherical one, because they are "non-intrinsic", i.e. lacking in absolute geometric meaning.

## 1.2 Drawbacks of Current Approaches

### 1.2.a The Infinite-Dimensional Approach

The standard mathematical presentation of the nonlinear filtering problem, as can be seen for example in Lipster and Shiryaev [12], and Pardoux [14], involves a measure-valued SDE called the Zakai equation (or the Fujisaki-Kallianpur-Kunita formula). This is virtually never used in real-time filtering applications because it is impossible to store enough data to update an infinite-dimensional SDE (although see Lototsky, Mikulevicius, and Rozovskii [13] for a computational method using a Wiener chaos expansion).

### 1.2.b Finite-Dimensional Filters

Under certain circumstances, the conditional law can be described using a finite set of parameters. Although this topic is outside the scope of this article, an account of recent progress using geometric methods can be found in Cohen de Lara [3]. Apart from the Kalman filter, these methods are not widely used in practice, since the parameters may be difficult to determine in theory, large in number, and difficult to update computationally.

### 1.2.c The Extended Kalman Filter and Second-order Filters.

Linearizing the state and observation about the most recent state estimate, and then applying the Kalman Filter, gives the Extended Kalman Filter; see Jazwinski [11] and Bar-Shalom and Fortmann [1]. The goal here is no longer to describe the full conditional distribution of the state given the observations, but merely to approximate the conditional expectation and the conditional covariance. As mentioned above, the output is coördinate-dependent. A careful asymptotic analysis of this and other approximation schemes has been given by Picard [15] - see also references therein.

### 1.3 Desirable Properties of a Nonlinear Filtering Algorithm

#### 1.3.a State Evolves Continuously, Observations are Discrete

The state dynamics (for example, the dynamics of an aircraft) should be modeled by a stochastic process $\{X_t, t \geq 0\}$ in **continuous time**, on a differentiable manifold $N$. However since digital implementation of a filtering algorithm is carried out using discrete-time observations, the filter should involve observations $Y_1, Y_2, \ldots$ collected at **discrete times** $t_1 < t_2 < \ldots$ on another manifold $M$.

#### 1.3.b State Estimates Should Not Be Coördinate-Dependent

Let $\{X_t^{(1)}, t \geq 0\}$ and $\{X_t^{(2)}, t \geq 0\}$ be representations of $\{X_t, t \geq 0\}$ in two coördinate systems, where $X_t^{(2)} = \phi(X_t^{(1)})$. Likewise let $Y_1^{(1)}, Y_2^{(1)}, \ldots$ and $Y_1^{(2)}, Y_2^{(2)}, \ldots$ be representations of $Y_1, Y_2, \ldots$ in two coördinate systems. We require that our state estimate of $X_{t_n}^{(2)}$, given $\{Y_1^{(2)}, \ldots, Y_n^{(2)}\}$, be the image under $\phi$ of our state estimate of $X_{t_n}^{(1)}$, given $\{Y_1^{(1)}, \ldots, Y_n^{(1)}\}$. Notice carefully that this criterion **excludes** use the conditional expectation $E\left[X_{t_n}^{(1)} \middle| Y_1^{(1)}, \ldots, Y_n^{(1)}\right]$ as the state estimate, because it does **not** have this kind of invariance. The replacement of conditional expectation by an "intrinsic location parameter" is the main theoretical contribution of this work.

#### 1.3.c Must Coincide with the Kalman Filter in the Linear Case

When $\{X_t, t \geq 0\}$ is a continuous-time Gaussian process, and $Y_n$ is a linear function of $X_{t_n}$ with additive Gaussian noise, our filtering algorithm must give the Kalman filter state estimates (which fully describe the conditional distribution of the state, given the observations, in such a context.)

#### 1.3.d Optimality up to Some Order

When the noise covariance of $X$, and the observation covariance are taken to be $O(\gamma^2)$, where $\delta \approx \gamma^2$ is the time interval between observations, we are seeking an algorithm which is optimal up to $O(\gamma^3)$, in a sense to be made precise later.

#### 1.3.e Stability

The important issue of stability will not be studied here. For results on the stability of the EKF, see Bossanne et al [2] and Deza et al [7].

## 2 The Nonlinear Model and its Induced Geometry

The geometric ideas in this section may be unfamiliar to filtering theorists, so we shall illustrate them with reference to the specific example of Section 1.1.

### 2.1 The State Process

Consider a (possibly degenerate) Markov diffusion process $\{X_t, 0 \leq t \leq \delta\}$ on $N \cong R^p$, written in local coördinates as

$$dX_t^i = b^i(X_t)\,dt + \sum_{j=1}^{p} \sigma_j^i(X_t)\,dW_t^j, \quad i = 1, 2, \ldots, p, \tag{1}$$

where $\sum b^i \frac{\partial}{\partial x_i}$ is a vector field on $R^p$, $\sigma(x) \equiv (\sigma_j^i(x)) \in L(R^p; T_x R^p)$, and $W$ is a Wiener process in $R^p$. We assume for simplicity that the coefficients $b^i$, $\sigma_j^i$ are $C^3$ with bounded first derivative.

## 2.2   Geometry Induced by the State Process

Such a $\sigma$ induces a $C^2$ semi-definite metric $\langle .|. \rangle$ on the cotangent bundle, which we call the **diffusion variance semi-definite metric,** by the formula

$$\langle dx^i | dx^k \rangle_x \equiv (\sigma \cdot \sigma(x))^{ik} \equiv \sum_{j=1}^{p} \sigma_j^i(x) \sigma_j^k(x) . \tag{2}$$

This semi-definite metric is actually intrinsic: changing coördinates for the diffusion will give a different matrix $(\sigma_j^i)$, but the same semi-definite metric. The $p \times p$ matrix $((\sigma \cdot \sigma)^{ij})$ defined above induces a linear transformation $\alpha(x): T_x^*N \to T_x N$, i.e. from the cotangent space to the tangent space at $x$, namely

$$\alpha(x)(dx^i) \equiv \sum (\sigma \cdot \sigma)^{ij} \partial/\partial x_j . \tag{3}$$

Let us make a **constant-rank assumption**, i.e. that there exists a rank $r$ vector bundle $E \to N$, a sub-bundle of the tangent bundle, such that $E_x = \text{range}(\sigma(x)) \subseteq T_x N$ for all $x \in N$. Darling [5] presents a global geometric construction of a **canonical sub-Riemannian connection** $\nabla^\circ$ for $\langle .|. \rangle$, with respect to a generalized inverse $g$, i.e. a vector bundle isomorphism $g: TN \to T^*N$ such that

$$\alpha(x) \bullet g(x) \bullet \alpha(x) = \alpha(x) . \tag{4}$$

In local coördinates, $g(x)$ is expressed by a Riemannian metric tensor $(g_{rs})$, such that if $\alpha^{ij} \equiv (\sigma \cdot \sigma)^{ij}$, then

$$\sum_{r,s} \alpha^{ir} g_{rs} \alpha^{sj} = \alpha^{ij} . \tag{5}$$

The local connector $\Gamma(x) \in L(T_x R^p \otimes T_x R^p; T_x R^p)$ for $\nabla^\circ$ can be written as:

$$2g(\Gamma(x)(u \otimes v)) \cdot w = D\langle g(v)|g(w)\rangle(u) + D\langle g(w)|g(u)\rangle(v) - D\langle g(u)|g(v)\rangle(w), \tag{6}$$

where $g(\Gamma(x)(u \otimes v))$ is a 1-form, acting on the tangent vector $w$. This formula coincides with the formula for the Levi-Civita connection in the case where $\langle .|. \rangle$ is non-degenerate; for more details, see Darling [5]. Our connection $\nabla^\circ$ gives rise to notions of geometry such as geodesics, parallelism, covariant differentiation, exponential map, and curvature, as explained in texts such as Darling [4], Sakai [16]. We assert:

*Axiom A:*   *The appropriate geometry for the state process is the one induced by the diffusion variance semi-definite metric.*

## 2.3 Intrinsic Description of the Process

The **intrinsic** version of (1) is to describe $X$ as a diffusion process on the manifold $N$ with generator

$$L \equiv \xi + \frac{1}{2}\Delta \tag{7}$$

where $\Delta$ is the (possibly degenerate) Laplace-Beltrami operator associated with the diffusion variance, and $\xi$ is a vector field, whose expressions in the local coördinate system $\{x^1, ..., x^p\}$ are as follows:

$$\Delta = \sum_{i,j} (\sigma \cdot \sigma)^{ij} \{D_{ij} - \sum_k \Gamma^k_{ij} D_k\} \, , \, \xi = \sum_k \{b^k + \frac{1}{2}\sum_{i,j} (\sigma \cdot \sigma)^{ij} \Gamma^k_{ij}\} D_k \, . \tag{8}$$

Note that $\sum \Gamma^k_{ij} (\sigma \cdot \sigma)^{ij}$ has been specified by (6).

## 2.4 Target Tracking Example

### 2.4.a State Process

The state $x$ consists of a column vector whose components $(p, v, a) \in R^3 \times R^3 \times R^3$ are respectively the location, velocity, and acceleration of the target in three-dimensional space. We model the acceleration as an Ornstein-Uhlenbeck process, with the constraint that acceleration must be perpendicular to velocity, i.e. $v \cdot a = 0$, or equivalently that $\|v\|^2$ is a constant. The $(v, a)$ components can be considered as taking values in the four-dimensional manifold $TS^2 \subset R^6$.

Thus within a Cartesian frame, the equations of motion take the nonlinear form:

$$\begin{bmatrix} dp \\ dv \\ da \end{bmatrix} = \begin{bmatrix} 0 & I & 0 \\ 0 & 0 & I \\ 0 & -\rho(x)I & -\lambda P(v) \end{bmatrix} \begin{bmatrix} p \\ v \\ a \end{bmatrix} dt + \begin{bmatrix} 0 \\ 0 \\ \gamma P(v) \, dW(t) \end{bmatrix}, \tag{9}$$

where the square matrix consists of nine $3 \times 3$ matrices, $\lambda$ is a positive time constant, $\gamma \equiv \sqrt{\lambda c_1}$ determines the noise intensity,

$$\rho(x) \equiv \|a\|^2 / \|v\|^2, \tag{10}$$

$$P(v) \equiv I - \frac{vv^T}{\|v\|^2} \in L(R^3; R^3), \tag{11}$$

and $W$ is a three-dimensional Wiener process. Note that $P(v)$ is precisely the projection onto the orthogonal complement of $v$ in $R^3$, and $\rho(x)$ has been chosen so that $d(v \cdot a) = 0$. (D'Souza et al [9] describe a procedure for estimating $\lambda$, but in our simulations we assign to it a predetermined value.) The constancy of $\|v\|^2$ implies that

$$DP(v) \zeta_v = \frac{-1}{\|v\|^2} \{\zeta_v v^T + v \zeta_v^T\} \, . \tag{12}$$

### 2.4.b  Geometry Induced by the State Process

The diffusion variance metric (2) is degenerate here; noting that $P^2 = P$, we find

$$\alpha \equiv \sigma \cdot \sigma \equiv \begin{bmatrix} 0 & 0 & 0 \\ 0 & 0 & 0 \\ 0 & 0 & \gamma^2 P(v) \end{bmatrix}, \tag{13}$$

where 0 denotes $0_{3\times 3}$. The rescaled Euclidean metric $g = \gamma^{-1} I_9$ on $R^9$ is a generalized inverse to $\alpha$ in the sense of (4), because $P^2 = P$. In Section 5 of [5] we show in more detail that the corresponding local connector $\Gamma(x)$ as in (6) is given by

$$\Gamma(x)(\zeta \otimes \varsigma) = \frac{S(\zeta \otimes \varsigma)v}{2\|v\|^2}, \quad S(\zeta \otimes \varsigma) \equiv \begin{bmatrix} 0 \\ \zeta_a \varsigma_a^T + \varsigma_a \zeta_a^T \\ -\zeta_v \varsigma_a^T - \varsigma_v \zeta_a^T \end{bmatrix}. \tag{14}$$

where a tangent vector $\zeta$ to $R^9$ is broken down into three 3-dimensional components $\zeta_p, \zeta_v, \zeta_a$. Note

$$\Gamma(x)(\sigma \cdot \sigma(x)) = \frac{\gamma^2}{\|v\|^2} \begin{bmatrix} 0 \\ P(v) \\ 0 \end{bmatrix} v = \begin{bmatrix} 0 \\ 0 \\ 0 \end{bmatrix}. \tag{15}$$

### 2.4.c  The Intrinsic Vector Field

It follows from (8), (9), and (15) that the formula for the intrinsic vector field $\xi$ is:

$$\xi(x) = \begin{bmatrix} v \\ a \\ -\rho(x)v - \lambda P(v)a \end{bmatrix}. \tag{16}$$

Differentiate under the assumptions $\|v\|^2$ is constant and $v \cdot a = 0$, to obtain

$$D\xi(x)(\zeta) = \begin{bmatrix} 0 & I & 0 \\ 0 & 0 & I \\ 0 & \lambda Q - \rho & -\lambda P - 2Q \end{bmatrix} \begin{bmatrix} \zeta_p \\ \zeta_v \\ \zeta_a \end{bmatrix}, \quad Q \equiv \frac{va^T}{\|v\|^2}. \tag{17}$$

In Section 5 of [5], we show that, if we write a symmetric tensor $\chi \in T_x N \otimes T_x N$ in $3 \times 3$ blocks as the matrix

$$\begin{bmatrix} \chi_{pp} & \chi_{pv} & \chi_{pa} \\ \chi_{vp} & \chi_{vv} & \chi_{va} \\ \chi_{ap} & \chi_{av} & \chi_{aa} \end{bmatrix},$$

where $\chi^T_{av} = \chi_{va}$, then

$$D^2\xi(x)(\chi) = \frac{-2}{\|v\|^2}\begin{bmatrix} 0 \\ 0 \\ \text{Tr}(\chi_{aa})v + (\chi_{va} + \chi_{av})a \end{bmatrix}. \tag{18}$$

### 2.4.d Curvature of the State Space

The curvature tensor is given by the formula (omitting $x$):

$$R(\varsigma,\eta)\zeta = D\Gamma(\eta)(\zeta\otimes\varsigma) - D\Gamma(\varsigma)(\zeta\otimes\eta) + \Gamma(\Gamma(\zeta\otimes\varsigma)\otimes\eta) - \Gamma(\Gamma(\zeta\otimes\eta)\otimes\varsigma), \tag{19}$$

and this may easily be computed from (14), noting that, for example, since $\|v\|^2$ is constant,

$$D\Gamma(\eta)(\zeta\otimes\varsigma) = S(\zeta\otimes\varsigma)\eta_v / (2\|v\|^2).$$

## 2.5 Covariance Tensor of a Random Variable in a Riemannian Manifold

We now introduce a local covariance concept, which we use for describing the uncertainty in the state estimates. Suppose $N$ is **any** manifold with a torsion-free connection, $\mu \in N$, and $\exp_\mu$ is the exponential map from $T_\mu N$ to $N$. Suppose $S$ is a random variable with values in $N$ (we assume that $S$ takes values in the image of the set on which $\exp_\mu$ is injective), and $\Sigma$ is a symmetric element of $T_\mu N \otimes T_\mu N$.

### 2.5.a Definition

*S will be said to **centered at $\mu$ with covariance tensor** $\Sigma$ if $\eta \equiv \exp_\mu^{-1} S$ satisfies $E[\eta] = 0 \in T_\mu N$, and for any cotangent vectors $\theta$ and $\lambda$ at $\mu$,*

$$E[(\theta\cdot\eta)(\lambda\cdot\eta)] = (\theta\otimes\lambda)\cdot\Sigma. \tag{20}$$

In more concrete terms, if $\{e_1, ..., e_p\}$ is some basis of $T_\mu N$, and

$$\Sigma = \sum_{i,j}\Sigma^{ij}e_i\otimes e_j,$$

then $(\Sigma^{ij})$ is the covariance matrix of the random vector $(\eta^1, ..., \eta^p)^T$ defined by $\exp_\mu^{-1} S = \sum_i \eta^i e_i$.

## 2.6 The Observation Covariance Metric

In our model, the observation $Y_n$ will be the image under the exponential map of a zero-mean random variable $V_n$ in the tangent space at $\psi(X_{t_n})$. Thus when $\psi(X_{t_n}) = y$, $Y_n$ is centered at $y$ with covariance tensor $\beta(y)$, a non-degenerate symmetric tensor in $T_y M \otimes T_y M$. Provided $y \to \beta(y)$ is sufficiently regular, it serves as the metric tensor for a metric $\langle.|.\rangle_o$ on the cotangent bundle of $M$, called the **observation covariance metric,** namely

$$\langle dy^i | dy^j \rangle_o = \beta^{ij}.$$

We assert:

***Axiom B:*** *The appropriate metric for the observation space is the observation covariance metric, not the Euclidean metric.*

We denote by $(g^o_{ij})$ the metric tensor on $TM$, inverse to $(\beta^{ij})$, and by $g^o$ the associated Riemannian metric. The Levi-Civita connection $\nabla$ on $M$ has local connector $\bar{\Gamma}(.)$, computed as follows:

$$\bar{\Gamma}^m_{ij} \equiv \frac{-1}{2}\sum_k \left\{ g^o_{jk}\frac{\partial \beta^{km}}{\partial y_i} + g^o_{ik}\frac{\partial \beta^{km}}{\partial y_j} + \frac{\partial g^o_{ij}}{\partial y_k}\beta^{mk} \right\} . \tag{21}$$

## 2.7   Target Tracking Example, Continued

### 2.7.a   Observation Function

The observables are respectively: range, angle from vertical, azimuth, and range-rate (all measured from a missile with known state $(p_M, v_M, a_M)$) and a fictitious measurement; the latter is a zero-mean Gaussian random variable representing a fictitious observation of the inner product of velocity and acceleration of the target, which according to our model should be zero. Take $\psi : R^3 \times R^3 \times R^3 \to R_+ \times S^2 \times R^2$ to be:

$$\psi(p, v, a) \equiv (\Phi(p - p_M), \|v - v_M\|, a \cdot v) \equiv ((y^1, y^2, y^3), y^4, y^5) \tag{22}$$

where $\Phi \equiv h^{-1}$, and $h$ is the spherical coördinate transformation

$$h(r, \theta, \phi) = (r\sin\theta\cos\phi, r\sin\theta\sin\phi, r\cos\theta) . \tag{23}$$

For the sake of brevity, we omit here the calculations of the first and second derivatives of $\psi$.

### 2.7.b   Observation Covariance Metric

The covariance matrix for the five observed quantities is taken to be of the form

$$\beta(y) \equiv \text{diag}\left( r^2 s_0, \frac{s_1}{r^2} + s_2, \frac{s_3}{r^2} + s_4, r^2 s_5, \sigma_F^2 \right) \equiv \text{diag}(h_1, h_2, h_3, h_4, h_5) , \tag{24}$$

where $r \equiv y^1$ denotes range, and other quantities are constants. Then

$$\frac{\partial \beta}{\partial y_i} = \delta_{1i}\text{diag}(h_1', h_2', h_3', h_4', 0) ,$$

and with $g^o(y) = \text{diag}(h_1^{-1}, h_2^{-1}, h_3^{-1}, h_4^{-1}, h_5^{-1})$,

$$\frac{\partial g^o}{\partial y_k} = -\delta_{1k}K, \ K \equiv \text{diag}\left( \frac{h_1'}{h_1^2}, \frac{h_2'}{h_2^2}, \frac{h_3'}{h_3^2}, \frac{h_4'}{h_4^2}, 0 \right).$$

Now (21) gives

$$\bar{\Gamma}^k_{ij} = \frac{1}{2}\left\{\delta_{ij}\delta_{1k}\frac{h_i'}{h_i^2}h_1 - \delta_{1i}\delta_{jk}\frac{h_k'}{h_j} - \delta_{1j}\delta_{ik}\frac{h_k'}{h_i}\right\}.$$

Let $f: R^5 \otimes R^5 \to R^5$ be the bilinear form with entries

$$f^k(u \otimes v) \equiv -\left(\frac{u^1 v^k + u^k v^1}{2}\right)\frac{h_k'}{h_k}, \quad k = 1, ..., 5.$$

The components of the connector are given by:

$$\bar{\Gamma}^k(u \otimes v) = f^k(u \otimes v) + \delta_{1k}\left(\frac{h_1}{2}\right)u^T K v, \quad k = 1, ..., 5. \tag{25}$$

## 2.8 Summary: the Model in Intrinsic Terms

Following the discussion above, we can rephrase the nonlinear filtering model in an abstract way.

### 2.8.a Model

The model consists of:

- A manifold $N$, called the state space, a canonical sub-Riemannian connection $\nabla^\circ$ induced by a diffusion variance semi-definite metric $\langle .|. \rangle$ on $T^*N$, and a vector field $\xi$ on $N$; these serve to define the generator $L \equiv \xi + \frac{1}{2}\Delta$ of a diffusion process $X$ on $N$;

- A Riemannian manifold $(M, g^o)$, called the observation space, and the Levi-Civita connection $\nabla$ induced by $g^o$.

- A $C^3$ function $\psi: N \to M$, called the observation function.

### 2.8.b Data

Data consist of:

- A point $\hat{\mu}_0 \in N$, called the initial state estimate;

- $\hat{\Sigma}_0 \in T_{\mu_0}N \otimes T_{\mu_0}N$, the covariance tensor of the initial state estimate;

- A sequence of times $0 < t_1 < t_2 < \ldots$, and for each $n \geq 1$ a noisy observation $Y_n$ of $\psi(X_{t_n})$ (in the sense of paragraph 2.6).

### 2.8.c Goal

The goal is to construct a sequence of state and covariance estimates $(\hat{\mu}_n, \hat{\Sigma}_n)$ for the state process, $n = 1, 2, \ldots$, with the following two properties:

- For a linear system subject to invertible smooth non-linear transformations, our estimates should be the transforms of the Kalman filter estimates.

- The construction of $(\hat{\mu}_n, \hat{\Sigma}_n)$ is intrinsic – i.e. unaffected by choice of coördinates – and optimal up to $O(\gamma^4)$, where $\gamma$ is a noise intensity parameter.

## 3    Intrinsic Geometric Quantities Associated with the Model

The following calculations will apply to every time interval $[t_{n-1}, t_n]$, but for notational simplicity we shall treat only use a time interval $[0, \delta]$. Let us fix a starting-point $x_0 \in N$, and $\Sigma_0 \in T_{\mu_0} N \otimes T_{\mu_0} N$.

### 3.1    Basic Notations

#### 3.1.a    The Deterministic Flow and Its Derivatives

$\{\phi_t, t \geq 0\}$ is the flow of the vector field $\xi$, and $x_t \equiv \phi_t(x_0)$. We assume non-explosion, so $\phi_t : N \to N$ is a $C^2$ diffeomorphism for each $t$. The flow $\{\phi_t, t \geq 0\}$ induces a two-parameter semigroup:

$$\tau_s^t \equiv d(\phi_t \bullet \phi_s^{-1})(x_s) \in L(T_{x_s} N; T_{x_t} N) . \tag{26}$$

In local coordinates, we compute $\tau_s^t$ as a $p \times p$ matrix, given by

$$\tau_s^t = \exp\{\int_s^t D\xi(x_u) \, du\} . \tag{27}$$

#### 3.1.b    Push-Forward of the Diffusion Variance Semi-definite Metric

For any vector field $\zeta$ on $N$, and any differentiable map $\phi : N \to P$ into a manifold $P$, the "push-forward" $\phi_* \zeta$ takes the value $d\phi \cdot \zeta(x) \in T_y P$ at $y \equiv \phi(x) \in P$; likewise $\phi_*(\zeta \otimes \zeta') \equiv \phi_* \zeta \otimes \phi_* \zeta'$.

The diffusion variance semi-definite metric $\langle . | . \rangle_{x_t}$ can be considered as an element of $T_{x_t} N \otimes T_{x_t} N$. Hence the following two quantities are intrinsic:

$$\Pi_t \equiv \Sigma_0 + \int_0^t (\phi_{-s})_* \langle . | . \rangle_{x_s} ds = \Sigma_0 + \int_0^t \tau_s^0 (\sigma \cdot \sigma)(x_s) (\tau_s^0)^T ds \in T_{x_0} N \otimes T_{x_0} N ; \tag{28}$$

$$\Xi_t \equiv (\phi_t)_* \Pi_t = \int_0^t \tau_s^t (\sigma \cdot \sigma)(x_s) (\tau_s^t)^T ds + \tau_0^t \Sigma_0 (\tau_0^t)^T \in T_{x_t} N \otimes T_{x_t} N . \tag{29}$$

#### 3.1.c    Second Fundamental Form

Given manifolds $N$ and $M$, with connections, whose local connectors are $\Gamma(.)$ and $\bar{\Gamma}(.)$, respectively, the second fundamental form $\nabla d\phi$ of a $C^2$ mapping $\phi : N \to M$ is a vector bundle morphism such that $\nabla d\phi(x) \in L(T_x N \otimes T_x N; T_{\phi(x)} M)$. It is expressed in local coordinates by:

$$\nabla d\phi(x)(v \otimes w) = D^2 \phi(x)(v \otimes w) - D\phi(x) \Gamma(x)(v \otimes w) + \bar{\Gamma}(y)(D\phi(x) v \otimes D\phi(x) w) \tag{30}$$

for $(v, w) \in T_x R^p \times T_x R^p$, $y \equiv \phi(x)$. We also need to compute $\nabla d\phi_\delta(x_0)$. From a formula in Darling [5],

$$\frac{\partial}{\partial t}(\tau_t^\delta D^2 \phi_t(x_0))(v \otimes w) = \tau_t^\delta D^2 \xi(x_t)(\tau_0^t v \otimes \tau_0^t w) ,$$

and together with (30), this leads to the expression:

$$\nabla d\phi_\delta(x_0)\,(v \otimes w) \;=\; \int_0^\delta \tau_t^\delta D^2\xi(x_t)\,(\tau_0^t v \otimes \tau_0^t w)\,dt - \tau_0^\delta \Gamma(x_0)\,(v \otimes w) + \Gamma(x_\delta)\,(\tau_0^\delta v \otimes \tau_0^\delta w)\,. \quad (31)$$

## 3.2  Approximate Intrinsic Location Parameter

Assume that $X$ is a diffusion on $N$ with generator $L \equiv \xi + \frac{1}{2}\Delta$, and random initial value $X_0$, centered at $x_0$ with covariance tensor $\Sigma_0$. A $C^3$ function $\psi : N \to M$ is given. We recall from [5] that there exists a vector in the tangent space $T_{\psi(x_\delta)}M$ which supplies a coördinate-independent replacement for the notion of expected value of $\psi(X_\delta)$. This vector, denoted $I_{x_0,\Sigma_0}[\psi(X_\delta)]$, is called the **approximate intrinsic location parameter (AILP)** of $\psi(X_\delta)$ in the tangent space $T_{\psi(x_\delta)}M$. We here omit any discussion of how the AILP is derived from the study of manifold-valued martingales, or its relation to harmonic mappings, but merely state the formula

$$I_{x_0,\Sigma_0}[\psi(X_\delta)] \;=\; \tfrac{1}{2}\{\nabla d\psi(x_\delta)\,(\Xi_\delta) + \psi_*\!\left[\nabla d\phi_\delta(x_0)\,(\Pi_\delta) - \int_0^\delta \tau_t^\delta \nabla d\phi_t(x_0)\,(d\Pi_t)\right]\}\,, \quad (32)$$

using the notations of Section 3.1. In the particular case where $\psi$ is the identity, we obtain $I_{x_0,\Sigma_0}[X_\delta]$, the AILP of $X_\delta$ in the tangent space $T_{x_0}N$, given by

$$I_{x_0,\Sigma_0}[X_\delta] \;=\; \tfrac{1}{2}\{\nabla d\phi_\delta(x_0)\,(\Pi_\delta) - \int_0^\delta \tau_t^\delta \nabla d\phi_t(x_0)\,(d\Pi_t)\}\,. \quad (33)$$

## 3.3  Numerical Evaluation of the Geometric Quantities Above

Suppose we have discretized the interval $[0,\delta]$. We now explain how to evaluate, at consecutive time steps $u < t$, the quantities $\{x_t, \tau_s^t, \Xi_t, m_t\}$, where $m_\delta \equiv I_{x_0,\Sigma_0}[X_\delta]$. The ODE $dx_t/dt = \xi(x_t)$ can be solved, for example, using the scheme

$$x_t \approx x_u + (t-u)\xi + \frac{(t-u)^2}{2}D\xi(x_u)\,(\xi) + \frac{(t-u)^3}{6}\{D^2\xi(x_u)\,(\xi \otimes \xi) + D\xi(x_u)\,(D\xi(x_u)\xi)\}\,, \quad (34)$$

where $\xi$ is short for $\xi(x_u)$. Using consecutive pairs $(x_u, x_t)$ computed from (34), we can discretize and solve (27), using the trapezium rule:

$$\tau_u^t \approx \exp\{\tfrac{t-u}{2}[D\xi(x_u) + D\xi(x_t)]\}\,, \quad (35)$$

where $\tau_0^0 = I$, and $\tau_0^t = \tau_u^t \tau_0^u$. Take $\Xi_0 = \Sigma_0$, and use (34), (35), and the trapezium rule to solve:

$$\Xi_t \approx \tfrac{t-u}{2}(\sigma \cdot \sigma)(x_t) + \tau_u^t\!\left[\Xi_u + \tfrac{t-u}{2}(\sigma \cdot \sigma)(x_u)\right](\tau_u^t)^T. \quad (36)$$

According to [5], the local formula for $m_\delta \equiv I_{x_0,\Sigma_0}[X_\delta]$ is

$$m_\delta \;=\; \tfrac{1}{2}\{\kappa_\delta - \tau_0^\delta \Gamma(x_0)\,(\Sigma_0) + \Gamma(x_\delta)\,(\Xi_\delta)\}\,, \quad (37)$$

$$\kappa_\delta \equiv \int_0^\delta \tau_t^\delta [D^2\xi(x_t)\,(\Xi_t) - \Gamma(x_t)\,(\sigma \cdot \sigma(x_t))]\,dt.$$

The last integral may be evaluated by:

$$\kappa_t \approx \frac{t-u}{2}[D^2\xi(x_t)\,(\Xi_t) - \Gamma(x_t)\,(\sigma \cdot \sigma(x_t))] + \tau_u^t\{\kappa_u + \frac{t-u}{2}[D^2\xi(x_u)\,(\Xi_u) - \Gamma(x_u)\,(\sigma \cdot \sigma(x_u))]\}$$

Finally (30), (32), and (33) show that

$$I_{x_0, \Sigma_0}[\psi(X_\delta)] = \frac{1}{2}\{D^2\psi(x_\delta)\,(\Xi_\delta) - J\Gamma(x_\delta)\,(\Xi_\delta) + \bar{\Gamma}(y_\delta)\,(J\Xi_\delta J^T)\} + Jm_\delta, \tag{38}$$

where $J \equiv D\psi(x_\delta)$. We compute $\nabla d\phi_\delta(x_0)$ in a similar way, using (31).

## 4   The Fundamental Theorems

Throughout this section, $X$ is a diffusion on $N$ with generator $L \equiv \xi + \frac{1}{2}\Delta$, and random initial value $X_0$, centered at $x_0$ with covariance tensor $\Sigma_0$. We are given a $C^3$ function $\psi: N \to M$, where $M$ is a Riemannian manifold of dimension $q$. Let $\beta(y) \in T_y M \otimes T_y M$ be the inverse metric tensor at $y \in M$, which can be interpreted as an observation covariance metric as in Section 2.6. Consider a single observation $Y_1$ of the form:

$$Y_1 \equiv \exp_{\psi(X_\delta)} V_1 \in M,$$

where $V_1$ is a mean-zero random vector in $T_{\psi(X_\delta)}M$, with covariance tensor $\beta(y)$ when $\psi(X_\delta) = y$, but which is otherwise independent of $U_0$ and of the Wiener process $W$.

### 4.1   Orders of Magnitude of Noise Terms

We shall suppose that, for some small parameter $\gamma$, the matrices for $\alpha \equiv \sigma \cdot \sigma$ (see (3)) and $\beta$ (see Section 2.6) satisfy

$$\alpha(x_t) = \gamma^2 \alpha_0(x_t),\ 0 \le t \le \delta;\ \beta(\psi(x_\delta)) = \gamma^2 \beta_0(\psi(x_\delta)); \tag{39}$$

where $\alpha_0$ is some other semi-definite metric, and $\beta_0$ another metric. Also assume that, with respect to the metric $g$ appearing in (4), the distribution of $U_0 \equiv \exp_{x_0}^{-1}(X_0)$ satisfies:

$$\|E[U_0]\| = O(\gamma^4),\ \|\Sigma_0\| \equiv \|\mathrm{Var}(U_0)\| = O(\gamma^2),\ \|E[T(U_0, U_0, U_0)]\| = O(\gamma^4), \tag{40}$$

for arbitrary tensor fields $T$ of type $(1, 3)$, whose norm is 1.

### 4.2   Definition

Let $U$ and $Z$ be integrable random variables in $R^p$, and $Y$ a random variable in $R^q$. We shall say that $Z$ approximates $E[U|Y]$ up to $O(\gamma^4)$ if

$$E[h(Y - E[Y]) \otimes (Z - U)] = O(\gamma^4), \tag{41}$$

for every $h \in C^1(R^q; R^p)$ with $\max\{\sup_y \|h(y)\|, \sup_y \|Dh(y)\|\} = 1$.

## 4.3 Interpretation of Tensors

To understand formulas such as (45) below, note that $\beta(y) \in T_y M \otimes T_y M$ can be interpreted as an element of $L(T_y^* M; T_y M)$, as in (3). For $y_\delta \equiv \psi(x_\delta)$, the adjoint of $J \equiv D\psi(x_\delta) \in L(T_{x_\delta} N; T_{y_\delta} M)$ is written $J^T \in L(T_{y_\delta}^* M; T_{x_\delta}^* N)$, and consequently $\Xi_\delta J^T \in L(T_{y_\delta}^* M; T_{x_\delta} N)$. For the convenience of users, (45) is expressed as a matrix product, but it actually represents a vector bundle morphism.

We quote the main result of [6].

## 4.4 Theorem (Intrinsic Conditional Expectation Formula)

*Consider the random vector* $U_\delta \oplus Z_\delta \in T_{x_\delta} N \oplus T_{\psi(x_\delta)} M$, *where* $y_\delta \equiv \psi(x_\delta)$, *given by*

$$U_\delta \equiv \exp_{x_\delta}^{-1}(X_\delta), \quad Z_\delta \equiv \exp_{y_\delta}^{-1}(Y_1). \tag{42}$$

*(i) Under the assumptions (39) and (40), the joint distribution of* $U_\delta$ *and* $Z_\delta$ *satisfies*

$$E\begin{bmatrix} U_\delta \\ Z_\delta \end{bmatrix} = \begin{bmatrix} I_{x_0, \Sigma_0}[X_\delta] \\ I_{x_0, \Sigma_0}[\psi(X_\delta)] \end{bmatrix} + O(\gamma^4); \tag{43}$$

*in terms of the approximate intrinsic location parameters of Section 3.2, where* $\Xi_\delta$ *is given by (29).*

*(ii)* $E[U_\delta | Z_\delta]$ *is approximated up to* $O(\gamma^4)$ *(in the sense of (41)) by*

$$I_{x_0, \Sigma_0}[X_\delta] + G\hat{Z}_\delta + \rho(\hat{Z}_\delta \otimes \hat{Z}_\delta) - E[\rho(\hat{Z}_\delta \otimes \hat{Z}_\delta)], \tag{44}$$

*where* $\hat{Z}_\delta \equiv Z_\delta - I_{x_0, \Sigma_0}[\psi(X_\delta)]$, *and* $G \in L(T_{y_\delta} M; T_{x_\delta} N)$ *is analogous to the Kalman gain, namely*

$$G \equiv \Xi_\delta J^T [J\Xi_\delta J^T + \beta(y_\delta)]^{-1}, \tag{45}$$

*where* $J \equiv D\psi(x_\delta)$, *and* $\rho \in L(T_{y_\delta} M \otimes T_{y_\delta} M; T_{x_\delta} N)$ *satisfies*

$$\rho(z \otimes z) \equiv \frac{1}{2}\{[I - GJ]\nabla d\phi_\delta(x_0)(\tau_\delta^0 Gz \otimes \tau_\delta^0 Gz) - G\nabla d\psi(x_\delta)(Gz \otimes Gz)\}, \tag{46}$$

$$E[\rho(\hat{Z}_\delta \otimes \hat{Z}_\delta)] = \rho(GJ\Xi_\delta) + O(\gamma^4).$$

*(iii)* $\text{Var}(U_\delta | Z_\delta)$ *is approximated up to* $O(\gamma^4)$ *by* $(I - GJ)\Xi_\delta$.

*(iv) If* $\hat{U}_\delta$ *denotes the difference between* $U_\delta$ *and (44), and if T is a tensor field of type* $(1, 3)$ *on N of norm 1, then*

$$\left\| E[T(\hat{U}_\delta, \hat{U}_\delta, \hat{U}_\delta)] \right\| = O(\gamma^4). \tag{47}$$

### 4.5    Computation of Exponential Barycenters

Let us recall from Emery and Mokobodzki [10] that an **exponential barycenter** for a random variable $S$ on a manifold $N$ with a torsion-free connection $\Gamma$ is a point $z \in N$ such that the random variable $\exp_z^{-1}(S) \in T_z N$ has mean zero. Suppose that are given a point $x \in N$, and moments

$$\mu \equiv E[\exp_x^{-1}(S)] \in T_x N, \quad \Sigma \equiv \text{Var}(\exp_x^{-1}(S)) \in T_x N \otimes T_x N. \tag{48}$$

We would like to compute from these moments an exponential barycenter for $S$, or at least a "good" approximation. We quote a result of Darling [6]. The norm $\|\cdot\|$ is with respect to some reference metric for $N$, which need not be related to the connection. Given the curvature tensor (19), the vector field $R\left(\frac{\partial}{\partial x_i}, \frac{\partial}{\partial x_j}\right)\frac{\partial}{\partial x_k}$ is denoted $R_{ijk}$.

### 4.6    Theorem (Exponential Barycenter Formula)

*Suppose that the moments (48) satisfy:* $\|\mu\| = O(\gamma)$, $\|\Sigma\| = O(\gamma^2)$, *for a small number $\gamma$. Define*

$$z \equiv \exp_x \{\mu - \frac{1}{3} \sum_{i,j,k} R_{ijk} \mu^i \Sigma^{jk}\}. \tag{49}$$

*Then* $E[\exp_z^{-1}(S)] = O(\gamma^4)$; *in other words, $z$ is an approximate exponential barycenter for $S$. If $T$ is a tensor field of type $(1, 3)$, and if* $\|E[T(\eta - \mu, \eta - \mu, \eta - \mu)]\| = O(\gamma^4)$, *where* $\eta \equiv \exp_x^{-1}(S) \in T_x N$, *then*

$$\left\| E[T(\exp_z^{-1}(S), \exp_z^{-1}(S), \exp_z^{-1}(S))] \right\| = O(\gamma^4). \tag{50}$$

### 4.7    Remark on the Validity of Recursion

Note carefully the relationship between the results (47) and (50), and the assumption (40). Assume that (40) holds. The conditional law of $X_\delta$, given $Y_1$, will be represented by a random variable $\hat{X}_\delta$ on $N$, whose exponential barycenter $z$ is computed according to (49), where $\mu$ is computed from (44), and $\Sigma$ is $(I - GJ)\Xi_\delta$. By (47) and (50), the random variable $\exp_z^{-1}(\hat{X}_\delta)$ in $T_z N$ will satisfy the same conditions that $U_0$ satisfied in (40). Therefore the algorithm can be repeated on the next time interval $[\delta, 2\delta]$, at the end of which we receive another observation $Y_2$, etc.

## 5    GI Filter Algorithm

As before, the state process $X$ is a diffusion on $N$ with generator $L \equiv \xi + \frac{1}{2}\Delta$, and random initial value $X_0$ centered at $\hat{\mu}_0$ with covariance tensor $\hat{\Sigma}_0$. A $C^3$ function $\psi : N \to M$ is given.

### 5.1    Discrete-time Observations

At each of the discrete times $0 < t_1 < t_2 < \ldots$, we make an observation $Y_n$ of the form:

$$Y_n \equiv \exp_{\psi(X_{t(n)})} V_n \in M,$$

where $V_n$ is a mean-zero random vector in $T_{\psi(X_{t(n)})}M$, with covariance tensor $\beta(y)$ when $\psi(X_{t(n)}) = y$, but which is otherwise independent of $X$. Define a sequence of sigma-fields

$$\mathfrak{I}_n^Y \equiv \sigma\{Y_1, ..., Y_n\}.$$

At time $t_n$, we wish to compute an $\mathfrak{I}_n^Y$-measurable random variable $\hat{\mu}_n$ with values in $N$, and an $\mathfrak{I}_n^Y$-measurable random variable $\hat{\Sigma}_n \in T_{\hat{\mu}_n}N \otimes T_{\hat{\mu}_n}N$, such that, conditional on $\mathfrak{I}_n^Y$, $X_{t_n}$ is approximately (i.e. up to $O(\gamma^4)$) centered at $\hat{\mu}_n$ with covariance tensor $\hat{\Sigma}_n$.

For any $n \geq 1$, suppose that $Y_1, ..., Y_{n-1}$ have been observed, from which we have calculated at time $t_{n-1}$ a state estimate $\hat{\mu}_{n-1}$ and its associated covariance tensor $\hat{\Sigma}_{n-1}$. The **GI Filter update formula** computes $\hat{\mu}_n$, and $\hat{\Sigma}_n$ as described in Sections 5.2 - 5.5.

### 5.2    Precomputation

First we carry out all the computations described in Section 3.3, starting from $x_0 \equiv \hat{\mu}_{n-1}$ and $\Sigma_0 \equiv \hat{\Sigma}_{n-1}$. The time interval $[t_{n-1}, t_n]$ is here represented as the time interval $[0, \delta]$, so when we refer to $x_\delta, \Xi_\delta$, etc., we are really referring to quantities at time $t_n$. The size of the computation depends on the number of sub-intervals into which we divide $[0, \delta]$, which can be as low as 1. Thus we obtain numerical expressions for all of the following:

$$x_\delta, J, y_\delta, \nabla d\psi(x_\delta), \Xi_\delta, \nabla d\phi_\delta(x_0), I_{x_0, \Sigma_0}[X_\delta], I_{x_0, \Sigma_0}[\psi(X_\delta)]$$

as well as $\tau_\delta^0 \equiv (\tau_0^\delta)^{-1}$. From these we compute the important coefficients $G$ and $\rho$, defined in (45) and (46), respectively.

### 5.3    Data Assimilation

All the formulas in this section are based on Theorem 4.4. We pull our new observation $Y_n$ back into the tangent space $T_{y_\delta}M$, by defining

$$Z_\delta \equiv \exp_{y_\delta}^{-1}(Y_n), \tag{51}$$

$$\hat{Z}_\delta \equiv Z_\delta - I_{x_0, \Sigma_0}[\psi(X_\delta)]. \tag{52}$$

See (58) for a simple formula for $Z_\delta$. In effect, $\hat{Z}_\delta$ is the "innovation", since it is the difference between the pulled-back observation $Z_\delta$ and its expected value, up to $O(\gamma^4)$. Next compute the approximate conditional expectation of $U_\delta \equiv \exp_{x_\delta}^{-1}(X_\delta)$, given $Y_n$, namely

$$\mu \equiv I_{x_0, \Sigma_0}[X_\delta] + G\hat{Z}_\delta + \rho(\hat{Z}_\delta \otimes \hat{Z}_\delta) - \rho(GJ\Xi_\delta) \in T_{x_\delta}N. \tag{53}$$

Note that $\rho$ is non-zero when any kind of non-linearity is present, so $\mu$ is a **quadratic** function of the innovation, not a linear one (as occurs in the Extended Kalman Filter, for example). The approximate conditional variance of $U_\delta \equiv \exp_{x_\delta}^{-1}(X_\delta)$, given $Y_n$, is

$$\Sigma \equiv (I - GJ) \, \Xi_\delta \in T_{x_\delta} N \otimes T_{x_\delta} N . \tag{54}$$

## 5.4    Update of the State Estimate

Recall that the canonical sub-Riemannian connection $\nabla^\circ$ on $N$ induces a curvature tensor, as in (19), and the vector field $R\left(\frac{\partial}{\partial x_i}, \frac{\partial}{\partial x_j}\right)\frac{\partial}{\partial x_k}$ is denoted $R_{ijk}$.

Define the new state estimate $\hat{\mu}_n$ by the formula of Theorem 4.6 for the conditional exponential barycenter of $X_\delta$, given $Y_n$, namely

$$\hat{\mu}_n \equiv \exp_{x_\delta}\{\mu - \frac{1}{3}\sum_{i,j,k} R_{ijk} \mu^i \Sigma^{jk}\} \in N . \tag{55}$$

See (59) for a simple formula for (55). Finally $\hat{\Sigma}_n \in T_{\hat{\mu}_n} N \otimes T_{\hat{\mu}_n} N$ is the push-forward of $\Sigma \in T_{x_\delta} N \otimes T_{x_\delta} N$ along the geodesic flow; see (60) for a straightforward method to compute $\hat{\Sigma}_n$.

## 5.5    Computation of Exponential Maps and Inverse Exponential Maps

Computation of $\exp_{x_\delta}(.)$ in (55), $\exp^{-1}_{y_\delta}(.)$ in (51), and $\hat{\Sigma}_n \in T_{\hat{\mu}_n} N \otimes T_{\hat{\mu}_n} N$, involves solving the first-order ordinary differential equations for the geodesic flows on the tangent bundles $TN$ and $TM$ respectively, which we describe briefly in Section 5.6. However there are also "single step" versions, which most practitioners will prefer to use, and which depend on the following classical formulas of local differential geometry, proved for example in Darling [6]:

### 5.5.a    Expansion of the Exponential Map

For $v \in T_x N \cong R^p$,

$$\exp_x v = x + v - \frac{1}{2}\Gamma(x)(v \otimes v) + \frac{1}{6}[2\Gamma(x)(\Gamma(x)(v \otimes v) \otimes v) - D\Gamma(x)(v)(v \otimes v)] + O(\|v\|^4) \tag{56}$$

### 5.5.b    Expansion of the Inverse Exponential Map

*The expansion for* $\exp^{-1}_y(z)$, *taking* $w \equiv z - y$, *is:*

$$w + \frac{1}{2}\Gamma(y)(w \otimes w) + \frac{1}{6}[D\Gamma(y)(w)(w \otimes w) + \Gamma(y)(\Gamma(y)(w \otimes w) \otimes w)] + O(\|w\|^4) . \tag{57}$$

### 5.5.c    Application

A simple way to approximate (51), avoiding use of the derivative of the connector, is:

$$Z_\delta \approx Y_n - y_\delta + \frac{1}{2}\bar{\Gamma}(y_\delta)((Y_n - y_\delta) \otimes (Y_n - y_\delta)) . \tag{58}$$

Since we have to differentiate $\Gamma(.)$ in any case to evaluate (19), take

$$v \equiv \mu - \frac{1}{3}\sum_{i,j,k} R_{ijk} \mu^i \Sigma^{jk} \tag{59}$$

and use (56) to compute (55). For the same choice of $v$, we may identify $T_v(T_{x_\delta}N)$ with $T_{x_\delta}N$, and so the derivative of $\exp_{x_\delta}(.)$ at $v$ can be viewed as a map from $T_{x_\delta}N$ to $T_{\hat{\mu}_n}N$ which is represented locally by a matrix $F$, where

$$F(w) \equiv w - \Gamma(x_\delta)(v \otimes w).$$

The formula for $F$ comes from differentiating (56) with respect to $v$. Finally the local formula for $\hat{\Sigma}_n$ is

$$\hat{\Sigma}_n \equiv F(I - GJ)\Xi_\delta F^T \in T_{\hat{\mu}_n}N \otimes T_{\hat{\mu}_n}N. \tag{60}$$

## 5.6 Geodesic Flow

This section is for those who seek more accurate calculations than the ones described in Section 5.5. The geodesic equation on $N$ can be represented as a first order ODE on the tangent bundle $TN$; in local coordinates, a solution is given by the geodesic flow

$$\begin{bmatrix} \gamma(s) \\ \zeta(s) \end{bmatrix} = \pi_s \left( \begin{bmatrix} \gamma(0) \\ \zeta(0) \end{bmatrix} \right) = \begin{bmatrix} \pi_s^H(\gamma(0), \zeta(0)) \\ \pi_s^V(\gamma(0), \zeta(0)) \end{bmatrix}, \quad 0 \leq s \leq 1, \tag{61}$$

in $R^p \oplus R^p$, satisfying the system of ODE

$$\begin{bmatrix} \gamma' \\ \zeta' \end{bmatrix} = h(\gamma, \zeta) \equiv \begin{bmatrix} \zeta \\ -\Gamma(\gamma)(\zeta \otimes \zeta) \end{bmatrix}. \tag{62}$$

(See Sakai [16], p. 56.) For example, to compute $\gamma(1) = \exp_{x_\delta}(v)$, the initial conditions will be $\gamma(0) = x_\delta$, $\zeta(0) = v$. Here $\{\gamma(s), 0 \leq s \leq 1\}$ will be a geodesic on $N$. In order to push a tensor forward from $T_{\gamma(0)}N$ to $T_{\gamma(1)}N$, i.e. to compute $(\pi_1^H)_*: T_{\gamma(0)}N \to T_{\gamma(1)}N$, we must compute the derivative flow $\{F(s), 0 \leq s \leq 1\}$ in $L(R^p \oplus R^p; R^p \oplus R^p)$ satisfying

$$F' = Dh(\gamma, \zeta)F, \quad F(0) = I. \tag{63}$$

We partition $F \equiv \begin{bmatrix} F_{11} & F_{12} \\ F_{21} & F_{22} \end{bmatrix}$ into $p \times p$ matrices; then $(\pi_1^H)_* = F_{11}(1)$.

## 6 Distinctive Features of the GI Filter

### 6.1 Invariance Under Change of Coördinate Systems

All the formulas for the GI Filter come from Theorems 4.4 and 4.6. All the mathematical quantities occurring in these two theorems are tensorial, i.e. are definable without using coördinates, and hence have the same intrinsic meaning for all coördinate systems. The only differences between computations in different coördinate systems will arise from numerical errors resulting from discretization, which can be made as small as desired.

## 6.2   Consistency with the Kalman Filter Under Nonlinear Transformations

Suppose that $\{X_t, t \geq 0\}$ is an Ornstein-Uhlenbeck process satisfying the SDE

$$dX_t = A_{n-1}X_t dt + \sigma_{n-1} dW_t, \quad t_{n-1} \leq t \leq t_n; \qquad (64)$$

here $A_{n-1}$ and $\sigma_{n-1}$ are $\mathfrak{I}^Y_{n-1}$-measurable $p \times p$ matrices. Also suppose that

$$Y_n = J_{n-1} X_{t_n} + V_n, \quad n = 1, 2, \ldots,$$

where $V_n \sim N_q(0, B_{n-1})$, $J_{n-1}$ and $B_{n-1}$ are $\mathfrak{I}^Y_{n-1}$-measurable matrices, and the $\{V_n\}$ are mutually independent random variables, independent of $W$. It is well known that the conditional distribution of $X_{t_n}$ given $\mathfrak{I}^Y_n$ is Gaussian, with conditional mean $\hat{\mu}^o_n$ and variance $\hat{\Sigma}^o_n$ given recursively by the Kalman Filter.

### 6.2.a   Proposition

*Suppose $\phi: R^p \to R^p$ and $\theta: R^q \to R^q$ are any $C^2$ diffeomorphisms. When the GI Filter is applied to the process $\{\phi(X_t), t \geq 0\}$, with observations $\theta(Y_1), \theta(Y_2), \ldots$ at times $t_1 < t_2 < \ldots$, the state estimator of $\phi(X_{t_n})$, given $\theta(Y_1), \ldots, \theta(Y_n)$, is $\phi(\hat{\mu}^o_n)$, with conditional covariance tensor*

$$D\phi(\hat{\mu}^o_n) \hat{\Sigma}^o_n (D\phi(\hat{\mu}^o_n))^T.$$

The theorem says, in effect, that when a nonlinear system is a transformed version of a linear system, then the GI Filter estimates are similarly transformed versions of the Kalman Filter estimates, as we desired in Section 1.3.c.

**Proof:** Since every step in the GI Filter is coördinate-independent, it suffices to prove the theorem when $\phi$ and $\theta$ are both the identity. When $\{X_t, t \geq 0\}$ satisfies (64), then (1) holds with $\sigma(x)$ not depending on $x$, and $b(x)$ is of the form $Ax$, where $A$ stands for $A_{n-1}$ when $t \in [t_{n-1}, t_n]$. The connector $\Gamma$ is zero on $N \equiv R^p$, so $\xi(x) = Ax$, and $D^2\xi(x)$ is zero.

Likewise since the covariance tensor is constant on observation space, the connector $\bar{\Gamma}$ is zero, and since $\psi(x) = J_{n-1}x$ is linear, we have $D^2\psi = 0$. We see that $I_{\hat{\mu}_{n-1}, \hat{\Sigma}_{n-1}}[X_\delta] = 0$ in (37) and $I_{\hat{\mu}_{n-1}, \hat{\Sigma}_{n-1}}[\psi(X_\delta)] = 0$ in (38). Abbreviating $J_{n-1}, B_{n-1}$ to $J, B$, etc., (45) becomes

$$G \equiv \Xi_\delta J^T [J \Xi_\delta J^T + B]^{-1},$$

and $\rho = 0$ in (46). In the constant-metric case, $\exp_x v = x + v$. From (53) and (55),

$$\hat{\mu}_n = x_\delta + G(Y_n - y_\delta), \quad \hat{\Sigma}_n = (I - GJ)\Xi_\delta,$$

which are the Kalman Filter estimates.   ◊

## 6.3   Relationship of the GI Filter with Standard Second-order Filters

Suppose that there is a coördinate system on state space in which the local connector $\Gamma(.)$, defined in (6), is zero, and a coördinate system on observation space in which $\bar{\Gamma}(.)$ (the Levi-Civ-ita connector for *M*) is zero. This situation occurs, for example, if $\sigma(x)$ does not depend on *x*, and $\beta(y)$ does not depend on *y*. The remaining nonlinearities come from $D^2\xi$, which need not be zero, and from $D^2\psi$. The formula (45) for *G* is comparable to standard formulas, but $\rho$ in (46), and the AILPs in (37) and (38) will be non-zero. The update formulas (53) and (55) become

$$\hat{\mu}_n \equiv x_\delta + I_{x_0, \Sigma_0}[X_\delta] + G\hat{Z}_\delta + \rho(\hat{Z}_\delta \otimes \hat{Z}_\delta) - \rho(GJ\Xi_\delta),$$

$$\hat{Z}_\delta \equiv Y_n - y_\delta - I_{x_0, \Sigma_0}[\psi(X_\delta)].$$

This quadratic formula is unlike the linear formulas found in the continuous-discrete Extended Kalman Filter (see [11], Theorem 8.1, p. 278), and the Truncated and Gaussian Second-order Filters ([11], equation (9.40), p. 345).

### 6.4     Size of Filtering Errors

We have said nothing about whether the observation function $\psi$ has properties (such as the rank of its derivative) sufficient to prevent filtering errors from diverging. See Picard [15] for a rigorous discussion of this point for a certain class of filters, under additional assumptions.

## 7     Software Implementation: Statistical Results

MATLAB codes for the GI filter and for a continuous-discrete Extended Kalman Filter have been developed for the tracking problem described in Sections 2.4 and 2.7. Although the GI Filter shows reduced bias, this example is not well suited to a statistical comparison, since results depend on many parameters, on the control law for the tracker, and on the coördinate systems chosen for the EKF. Moreover in a nine-dimensional state space, computations are relatively slow for both filters.

Consequently a much simpler example was selected as a context for statistical comparison. Here both the state process and the observations are one-dimensional, and the noise variance and observation variance do not vary, which forces both the local connectors $\Gamma(.)$ and $\bar{\Gamma}(.)$ to be zero. The model is:

$$dX_t = \xi(X_t)\,dt + \sqrt{\alpha}\,dW_t, \quad t \geq 0, \tag{65}$$

$$Y_n = \psi(X_n) + \sqrt{\beta}\,V_n, \quad n = 1, 2, \ldots, \tag{66}$$

where $\xi(x) \equiv -x^3/2$ and $\psi(x) \equiv x/(p + x^2)$. If one were using such a simple model in real life, direct calculation of the conditional density would be a natural approach in practice. The GI Filter and EKF were programmed merely for the sake of statistical comparison.

Note that $\psi$ has critical points at $\pm\sqrt{p}$, at which any filter is bound to perform badly. We chose parameters $p \equiv 0.1$, $\alpha \equiv 0.01$, and $\beta \equiv 0.001$, which cause the state process, which is a positive recurrent diffusion, to visit the critical points fairly often. Since the model is stationary, statistical characteristics of a filter may be observed by simply running the filter over thousands of cycles.

Fortunately all the formulas defined in Section 3 can be computed analytically, without recourse to numerical integration. For example, $x_t = x_0 / \sqrt{1 + x_0^2 t}$, and the AILP of $X_\delta$ is given by

$$m_\delta = \frac{-3}{2} \left\{ \frac{\alpha}{12 x_\delta^3} + x_0^{-4} \left( \Sigma_0 - \frac{\alpha}{3 x_0^2} \right) x_\delta^3 - x_0^{-6} \left( \Sigma_0 - \frac{\alpha}{4 x_0^2} \right) x_\delta^5 \right\}.$$

Histograms for the absolute value of the filter error over 10,000 filter cycles are shown. To preserve stability of the GI Filter, we placed a "collar" over the second order term $\rho \, (\hat{Z}_\delta \otimes \hat{Z}_\delta - GJ\Xi_\delta)$ in (53) so that its magnitude could not exceed that of the first order term. Note that large filtering errors occur less frequently for the GIF than for the EKF. This reflects the fact that, when the GI Filter is "thrown off" by the nonlinearity of the observation function $\psi$, it recovers more quickly than the EKF does. An example is shown in the following time series of 100 filter cycles, using the same $X$ series.

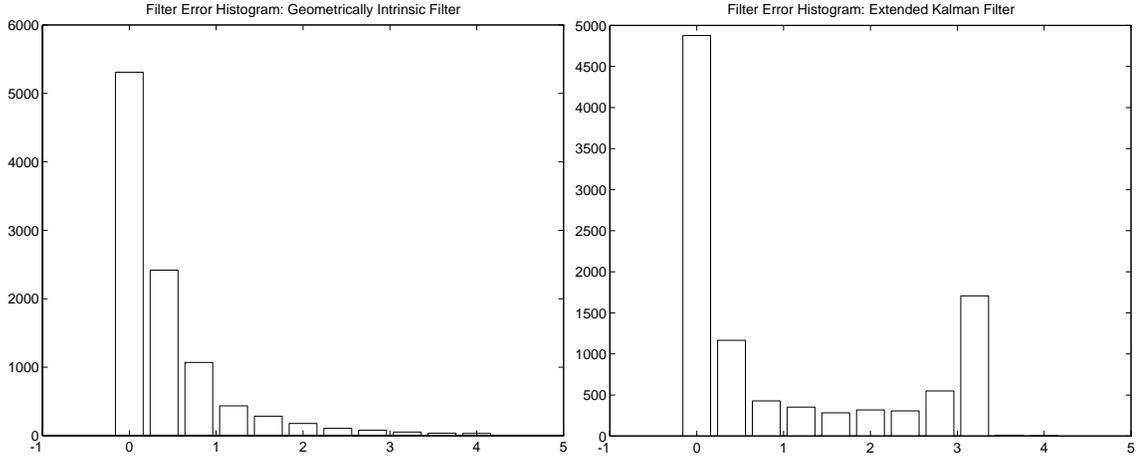

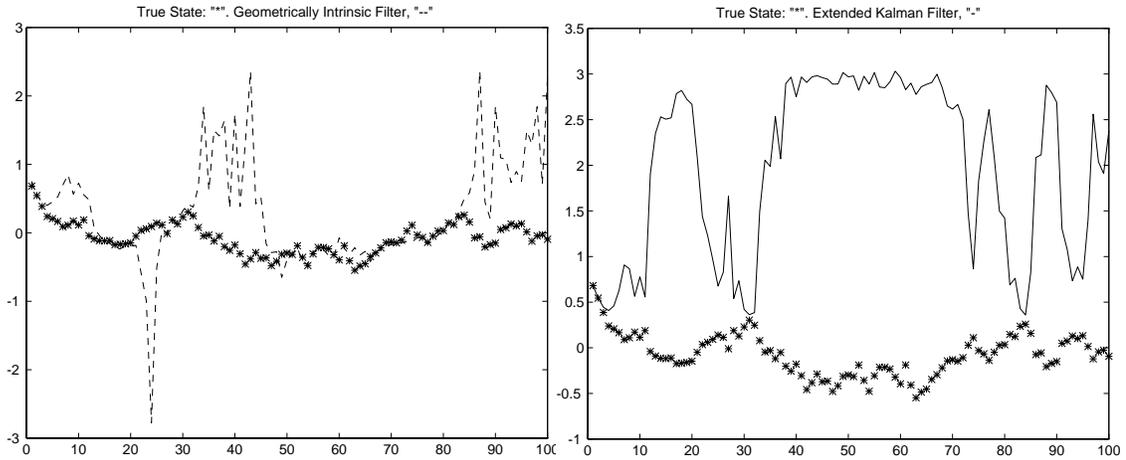

It should be emphasized that the parameter values chosen represent an extreme regime at the very edge of instability, and that for a large range of parameter values, for example for $p > 4$ and small α and β, the GI Filter and the EKF perform about the same. Moreover when the second order term $\rho\,(\hat{Z}_\delta \otimes \hat{Z}_\delta - GJ\Xi_\delta)$ is deleted, the performance advantage of the GI Filter disappears.

## 8 Conclusion

This project has demonstrated that there is a natural *intrinsic* generalization of the Kalman Filter to the fully nonlinear context, and that it is possible to implement it computationally. Unlike the Extended Kalman Filter, the state estimate for the GI Filter is a quadratic, not a linear, function of the observations. Computational experiments show that, when the observation function is highly nonlinear, there exist choices of the noise parameters at which the GI Filter significantly outperforms the EKF.

**Acknowledgment:** The author thanks Dr. James R. Cloutier of the Wright Laboratory for explaining the tracking problem analyzed above, and for valuable advice and encouragement. He also thanks Boris Rosovskii, Etienne Pardoux, and Ofer Zeitouni for helpful conversations, and the Statistics Department, University of California at Berkeley, for its hospitality during the writing of this article.